\documentclass[12pt,onecolumn,draftcls]{IEEEtran}      

\usepackage{tipa}

\usepackage{graphics}	
\usepackage{amsmath}	
\usepackage{amssymb}	
\usepackage{pifont}		
\usepackage{mathrsfs} 
\usepackage{mathtools}
\usepackage[dvipsnames]{xcolor}
\usepackage{tikz}
\usepackage{cite}
\usepackage{breqn}
\usepackage{enumerate}
\PassOptionsToPackage{hyphens}{url}
\usepackage[	
  		colorlinks = true, 
  		allcolors=black, 
		urlcolor=black,
	   ]{hyperref} 
\hypersetup{breaklinks=true}
\usepackage{url}

\usepackage[normalem]{ulem}

\newcommand\soutpars[1]{\let\helpcmd\sout\parhelp#1\par\relax\relax}
\long\def\parhelp#1\par#2\relax{%
  \helpcmd{#1}\ifx\relax#2\else\par\parhelp#2\relax\fi%
}

\def\tr{\mathrm{tr}}
\def\de{\mathrm{d}}
\def\R{\mathbb{R}}
\def\C{\mathbb{C}}
\def\D{\mathbb{D}}
\def\Z{\mathbb{Z}}
\def\T{\mathbb{T}}

\def\L{\mathcal{L}}
\def\B{\mathcal{B}}
\def\H{\mathcal{H}}

\def\S{\mathcal{S}}
\def\K{\mathcal{K}}

\newcommand{\diag}{\mathrm{diag}}

\newtheorem{proposition}{Proposition}
\newtheorem{theorem}{Theorem}

\newtheorem{definition}{Definition}
\newtheorem{remark}{Remark}

\newtheorem{example}{Example}

\definecolor{Royalblue}{cmyk}{1,0.30,0.2,0.2}

\title{
	Conal Distances Between \\ Rational Spectral Densities
}

\author{ Giacomo Baggio, Augusto Ferrante, and Rodolphe Sepulchre
\thanks{G.~Baggio and A.~Ferrante are with the Dipartimento di Ingegneria dell'€™Informazione, Universit\`a di Padova, via Gradenigo, 6/B€" I-35131 Padova, Italy. E-mail: \texttt{giacomo.baggio@studenti.unipd.it}, \texttt{augusto@dei.unipd.it}. R.~Sepulchre is with the Department of Engineering, University of Cambridge, Trumpington Street, Cambridge CB2 1PZ, UK. E-mail: \texttt{r.sepulchre@eng.cam.ac.uk}.
The research leading to these results has received funding from the European Research Council under the
Advanced ERC Grant Agreement Switchlet n.670645.}
}
\allowdisplaybreaks

\begin{document}

\maketitle
\thispagestyle{empty}
\pagestyle{empty}

\begin{abstract}
The paper generalizes Thompson and Hilbert metric to the space of spectral densities. The resulting complete metric space has the differentiable structure of a Finsler manifold with explicit geodesics. The resulting distances are filtering invariant, can be computed efficiently, and admit geodesic paths that preserve rationality; these are properties of fundamental importance in many engineering applications.
\end{abstract}


\section{Introduction}

In recent years the design and analysis of distances between spectral densities have received a renewed interest in the control and signal processing community (see \cite{B13,AV14} for two recent surveys on this topic). This interest primarily stems from a large number of applications in which the problem of quantifying dissimilarities between spectral densities is of crucial importance, such as spectral estimation \cite{GL03,FPR08,FMP12,EK08,E09,BETA,ALPHA,OPTIMAL_PREDICTION_ZORZI_2014}, speech processing  \cite{GBGM80,GKK95,JTG08,JLG08,JKG08,JLG12}, and time-series clustering \cite{L05,BDDBM05,BDCDM06,CP08,ACRV12,LDM17}, to cite a few.

The design of distances with the aim of solving computational engineering problems is a rich topic because of the interplay between  mathematical, modelling,  and computational considerations. Mathematical considerations include endowing the underlying space with a differential manifold structure such that the distance between two points corresponds to the length of a minimal geodesic.  This is especially relevant when dealing with problems involving approximation, smoothing, and averaging of spectral densities, e.g., in the context of speech morphing \cite{JLG08,JTG08,JLG12}.  The classical framework is Riemannian geometry, in which the  differential structure involves an inner product. The present paper uses the broader framework of Finsler geometry, where the differential structure only requires a norm. Modelling considerations include endowing the distance with suitable invariance properties, such that the mathematical distance is consistent with what is modelled. In the context of spectral densities,  filtering invariance emerges as a property that should hold when spectral densities model second-order stationary stochastic processes.  Invariance properties are receiving increasing attention in engineering, because they tend to make algorithms less sensitive to modelling assumptions.  Computational considerations include the existence of an algorithmic framework to perform the calculations necessary to the considered engineering problem, starting with the evaluation of the distance itself. They are of primary importance in high dimensional problems, and, a fortiori, so for infinite dimensional objects such as spectral densities. 

The starting point in this paper is to acknowledge that the space of spectral densities is a \emph{cone} and to revisit two classical distances that have been studied in cones: the part metric (often called  Thompson metric) and the projective metric introduced by Hilbert.   Applying this classical framework to the space of {\it rational} spectral densities, which seems novel, we show that the resulting metrics have a number of particularly desirable properties: 
\begin{itemize}
\item they are filtering invariant, a natural and desirable property in many applications; 
\item their calculations boil down to evaluating the $\H_{\infty}$-norm of minimum-phase spectral factors, a problem extensively studied in the control literature that can be performed via efficient routines;
\item they endow the cone of spectral densities with a Finsler geometry featuring explicit geodesic paths that can be chosen to be \emph{rational}.
\end{itemize}
In particular, we show that the Thompson metric is a close relative of a Riemannian metric recently studied in  \cite{JNG12}. However, the latter does not enjoy all the above-listed properties.

\emph{Paper structure.}  After some preliminary definitions and necessary notations, in Section~\ref{sec:motivations}, we illustrate some applicative scenarios that motivate the introduction of the new metrics. Section~\ref{sec:dist-cones} reviews the Finsler geometry of cones, with a special emphasis on the cone of positive definite matrices. Section~\ref{eq:Finslerian-dist} applies this geometry to the cone of rational spectral densities and discuss the properties inherited by this geometry.
Section~\ref{sec:appl} presents an applicative example of spectral interpolation via geodesic paths of the introduced metric, in the context of speech processing. Lastly, Section~\ref{sec:conclusions} collects some concluding remarks and future research directions.

\emph{Notation.} As usual, we denote by $\R$, $\C$, $\R^{n\times m}$, $\C^{n\times m}$, and $\R^{n\times m}(z)$, the set of real numbers, complex numbers, $n\times m$ real matrices, $n\times m$ complex matrices, and $n\times m$ real matrix-valued rational functions in $z\in\C$, respectively. $\T$ and $\D$ will denote, respectively, the unit circle, and the open unit disk in the complex plane. $\R_*^{n\times n}(z)$ will denote the set of $n\times n$ real matrix-valued rational functions of full rank on $\T$; it forms a multiplicative group. Given $A\in\C^{n\times n}$, $A^{\top}$, $A^{*}$, $\tr(A)$, and $\|A\|_F$ stand for the transpose, Hermitian transpose, trace and Frobenius norm of $A$, respectively. Let $\L^{n}_{2}[-\pi,\pi]$ be the space of $n$-dimensional vector-valued functions on $\T$ that are square integrable w.r.t. the normalized Lebesgue measure. The space $\L^{n}_{2}[-\pi,\pi]$ endowed with the inner product $\langle f,g\rangle_{2}:=\int_{-\pi}^{\pi}f(e^{j\theta})^{*}g(e^{j\theta})\,\frac{\de \theta}{2\pi}$, $f,g\in\L^{n}_{2}[-\pi,\pi]$, forms an Hilbert space. We let  $\S_{+}^{n\times n}$ denote the cone of $n\times n$ positive semi-definite Hermitian matrices and $\S_{+}^{n\times n}(\T)$ the cone of $n\times n$ bounded positive self-adjoint operators on $\L^{n}_{2}[-\pi,\pi]$, namely,
\begin{align*}
	\S_{+}^{n\times n}(\T):=\{\,&\Phi\colon \T\to\C^{n\times n} : \Phi(e^{j\theta})=\Phi(e^{j\theta})^{*}, \ \forall e^{j\theta}\in\T,  \\ &\ \ \ \ \   \text{ and } \langle f,\Phi f \rangle_{2}\ge 0, \ \forall f\in\L^{n}_{2}[-\pi,\pi] \,\}.
\end{align*}
Henceforth, elements of $\S_{+}^{n\times n}(\T)$ will be thought of as discrete-time \emph{spectral densities}\footnote{We remark that the standard definition of spectral density needs integrability only, see, e.g., \cite{LP15}. In this paper, however, we restrict the attention to a subset of ``well-behaved'' spectral densities that includes the (particularly important) class of bounded rational spectral densities.} and $\S_{+,\mathrm{rat}}^{n\times n}(\T)\subset \S_{+}^{n\times n}(\T)$ will denote the subset of $n\times n$ real \emph{rational}  bounded discrete-time spectral densities. 
Given $A(z)\in\R^{n\times m}(z)$, we let $A^{*}(z):=A^{\top}(1/z)$ and we use $A^{-*}(z)$ as a shorthand for $[A^{*}(z)]^{-1}$. A rational matrix $W(z)\in\R^{n\times r}(z)$ is called a spectral factor of the spectral density $\Phi(z)\in\S_{+,\mathrm{rat}}^{n\times n}(\T)$, if it satisfies $\Phi(z)=W(z)W^{*}(z)$. If the spectral factor $W(z)$ is analytic in (an open set containing) the complement of $\D$ with (generalized) inverse analytic in the complement of the closure of $\D$, then $W(z)$ is called minimum-phase (or canonical) spectral factor. The minimum-phase spectral factor of a rational spectral density $\Phi(z)\in\S_{+,\mathrm{rat}}^{n\times n}(\T)$ always exists and is essentially unique, that is, is unique up to post-multiplication by orthogonal $r\times r$ matrices, where $r$ is the normal rank (that is, the a.~e.~rank on $\mathbb{T}$) of $\Phi$ \cite{baggio2016factorization}. To conclude, we recall that, given a $n\times m$ matrix-valued function $G\colon\T\to\C^{n\times m}$ that is (essentially) bounded on $\T$, the $\L_\infty$-norm of $G$ is defined as $\|G\|_{\L_\infty}:=\mathrm{ess}\sup_{\theta\in[-\pi,\pi]}\sigma_{\max}(G(e^{j\theta}))$ where $\sigma_{\max}(A)$ denotes the largest singular value of $A\in\C^{n\times m}$. If, furthermore, $G$ can be analytically extended in an open set that contains the complement of $\D$, then its $\L_\infty$-norm coincides with the $\H_\infty$-norm $\|G\|_{\H_\infty}:=\sup_{z\in \C\,:\, |z|> 1}\sigma_{\max}(G(z))$.

\section{Motivational examples}\label{sec:motivations}

As mentioned in the introduction, engineering applications motivate  the design of distances that  (i) are computable, (ii) possess a differential structure amenable to interpolation and extrapolation, and (iii) are invariant to pre-processing of the data. In what follows, we  illustrate those concrete motivations in  representative specific applications.
\begin{example}[Time-series clustering \cite{L05,BDDBM05,BDCDM06}]
Let $\{y_{i}(k)\}_{k\in M_{i}}$, $M_{i}\subset\Z$, $i=1,2,\dots,N$, be a set of $p$-dimensional time series data collected from measurements and representing the (noisy) behavior of some physical or engineering systems. Each time series is modelled by a dynamical model estimated from the raw data. The estimated model is usually taken to be linear and time-invariant, and admits the (formal) input-output representation
\begin{equation}\label{eq:lin-mod-cluster}
	y_{i}(k) = W_{i}(z) e_{i}(k), \quad i=1,2,\dots,N,
\end{equation}
where $W_{i}(z)$ is a $p\times m$ rational matrix and $\{e_{i}(k)\}_{k\in\Z}$ an $m$-dimensional white noise process.
In many practical applications ranging from econometrics to biology, one needs to classify the different times series data $\{y_{i}(k)\}_{k\in M_{i}}$ in different groups, based on their similarity. This procedure is commonly known as \emph{time-series clustering}. 

The desired clustering can be performed by defining a suitable distance index between corresponding linear models in \eqref{eq:lin-mod-cluster}. Since the spectral density of the $i$-th time-series is $\Phi_{i}(z)=W_{i}(z)W_{i}^{*}(z)$, the distance can be defined in the space of (rational) spectral densities. In many real scenarios, the amount of time series data is massive and one seeks a distance that is, at the same time, computationally tractable and accurate.

\end{example}

\begin{example}[``Geodesic'' speech morphing\cite{JTG08,JLG12}] Consider two digital speech sources coming from two different individuals, say A and B, and described by time series $\{y_{\text{A}}(k)\}_{k=1}^{N}$ and $\{y_{\text{B}}(k)\}_{i=1}^{N}$, respectively. The \emph{morphing} of these two audio signal consists of gradually deforming the speech signal of A into the speech signal of B, creating a new hybrid signal that should preserve the speech-like quality and content of the original sources. Nowadays, there are many applications that benefit from such algorithms, especially in the areas of multimedia engineering and entertainment. 

By segmenting the  two time series into $M<N$ approximately stationary fragments, one can first estimate the spectral density of each fragment, say $\{\Phi_{\text{A},i}\}_{i=1}^{M}$ and $\{\Phi_{\text{B},i}\}_{i=1}^{M}$, and then perform morphing using a \emph{path} connecting the spectral densities of each fragment. A {\em geodesic} between $\Phi_{\text{A},i}$ and $\Phi_{\text{B},i}$ is a particularly convenient and natural choice for such a path. Of course, this requires the definition of a suitable metric in the space of spectral densities that equips the latter space with a differential structure. Thus, a geodesic interpolation path  between $\Phi_{\text{A},i}$ and $\Phi_{\text{B},i}$ provides a geometric solution to the morphing engineering question.
\end{example}

\begin{example}[THREE-like spectral estimation \cite{GL03,FPR08,FMP12,EK08,E09}] \label{ex:THREE-est}
Let $\{y(t)\}_{t\in\Z}$ be an $n$-dimensional zero-mean stationary process and let $\Psi\in\S_{+}^{m\times m}(\T)$ be an a priori estimate of the unknown spectral density of this process. Consider the bank of filters described by the transfer function
\[
	G(z)=(zI-A)^{-1}B, \quad A\in\R^{n\times n}, B\in\R^{n\times m},
\]
with $A$ strictly Schur stable, $B$ of full column rank, and the pair $(A,B)$ reachable. We assume to have an estimate of the asymptotic state covariance $\Sigma \in \S_{+}^{n\times n}$ of the system with transfer function $G(z)$ and input the process $\{y(t)\}_{t\in\Z}$.  The task is to estimate the spectral density of $\{y(t)\}_{t\in\Z}$ based on the available information.  Typically, the prior $\Psi$ is not consistent with the state covariance $\Sigma$. Therefore, it is necessary to find a spectral density which is as closest as possible, in some suitable sense, to $\Psi$, and, additionally, satisfies the ``consistency'' condition
\[
	 \int_{-\pi}^{\pi} G(e^{j\theta})\Phi(e^{j\theta}) G^{*}(e^{j\theta})\frac{\de \theta}{2\pi}=\Sigma.
\]
This formulation leads to  the following constrained optimization problem
\begin{align*}
\begin{split}
& \, \min_{\Phi\in \S_{+}^{m\times m}(\T)} d(\Psi,\Phi)\\
&\text{s.t. } \int_{-\pi}^{\pi} G(e^{j\theta})\Phi(e^{j\theta}) G^{*}(e^{j\theta})\frac{\de \theta}{2\pi}=\Sigma
\end{split}
\end{align*}
where $d\colon \S_{+}^{m\times m}(\T)\times \S_{+}^{m\times m}(\T)\to [0,\infty)$ is a suitable distance function in the cone of spectral densities.

One crucial aspect in the above estimation problem concerns the choice of the distance measure $d(\cdot,\cdot)$ to minimize. In order to be effective, this distance should satisfy some properties that naturally arise from the estimation setting. For instance, one natural requirement is that the distance between the prior process (described by spectral density $\Psi$) and the unknown process must be left unchanged if we filter both processes via the same filter. Another often desired property is that the distance is projective meaning that it is unaffected by scalings (in this case the ``shape'' of the unknown spectrum is actually estimated). Finally, the distance must be amenable to algorithmic optimization, which benefits from properties such as convexity and requires efficient numerical estimation of the distance and its first or second derivatives.


\end{example}

\section{Distances in cones}\label{sec:dist-cones}

Let $\K$ be a closed, solid, pointed, convex cone defined in a real Banach space $\B$ with norm $\|\cdot\|_{\B}$, that is, a closed subset $\K$ with the properties that: (i) the interior of $\K$, denoted by $\mathring{\K}$, is non-empty, (ii) $\K+ \K \subseteq \K$, (iii) $\K\cap-\K = \{0\}$, (iv) $\lambda\K\subseteq \K$ for all $\lambda\geq 0$. The cone $\K$ induces a partial ordering $\le_{\K}$ on $\B$ by
\[
	x\leq_{\K} y \iff y-x\in\K.
\]
For $x,y\in\K$, we say that $y$ dominates $x$ if there exists $\beta>0$ such that $x\le_{\K} \beta y$. We write $x\sim_{\K} y$ if $y$ dominates $x$, and $x$ dominates $y$. The relation $\sim_{\K}$ is an equivalence relation on $\K$. The corresponding equivalence classes are called \emph{parts} or \emph{components} of $\K$.

Given two elements $x,y\in \K\setminus \{0\}$, we define the following quantities
\begin{align}\label{eq:M}
	M(x,y)&:=\inf\{\lambda\,:\, x\le_{\K} \lambda y \},
\end{align}
if the set is non-empty, and $M(x,y):=\infty$ otherwise, and
\begin{align}\label{eq:m}
	m(x,y)&:=\sup\{\mu\,:\, \mu y\le_{\K} x\}=\frac{1}{M(y,x)}.
\end{align}

\begin{definition}[\hspace{-0.005cm}\cite{B73,T63}]\label{def1}
The \emph{Hilbert (projective) metric} and the \emph{Thompson (part) metric} between elements $x,y\in \K\setminus \{0\}$ are defined respectively by
\begin{align}
	d_{H}(x,y) &:= \log \frac{M(x,y)}{m(x,y)},\label{eq:hilbert}\\
	d_{T}(x,y) &:= \log \max\left\{M(x,y),M(y,x)\right\},\label{eq:thompson}
\end{align}
if $x\sim_{\K} y$, and $d_{H}(x,y)=d_{T}(x,y):=\infty$, otherwise.
\end{definition}

As a simple example, consider $\B=\R^{n}$ and $\K$ to be the positive orthant of $\R^{n}$, i.e. $\K:=\{(x_{1},\dots,x_{n})\,:\, x_{i}\geq 0,\, 1\le i\le n\}$. In this case, for $x,y\in\mathring{\K}$, it holds
\begin{align*}
	M(x,y) &= \max_{i}\{x_{i}/y_{i}\},\\
	m(x,y) &= \min_{i}\{x_{i}/y_{i}\},
\end{align*}
so that Hilbert and Thompson metrics on $\mathring{\K}$ read, respectively, as
\begin{align*}
	d_{H}(x,y) &= \log\frac{\max_{i}\{x_{i}/y_{i}\}}{\min_{i}\{x_{i}/y_{i}\}},\\
	d_{T}(x,y) &= \log\max\left\{\max_{i}\{x_{i}/y_{i}\},\max_{i}\{y_{i}/x_{i}\}\right\}.
\end{align*}

Thompson metric is a bona fide distance\footnote{We recall that a bona fide metric or distance function on a set $X$ is a function $d\colon X\times X \to [0,\infty)$ satisfying the following conditions for all $x,y,z\in X$: (i) $d(x,y)\ge 0$, (ii) $d(x,y)=0$ if and only if $x=y$, (iii) $d(x,y)=d(y,x)$, (iv) $d(x,z)\le d(x,y)+d(y,z)$.} on each part of the cone $\K$ (and, in particular, on the interior $\mathring{\K}$). Each part of $\K$ is a \emph{complete} metric space with respect to this metric provided that $\K$ is normal, i.e., there exists $\gamma>0$ such that $\|x\|_{\B}\le \gamma\|y\|_{\B}$ holds whenever $0\le_{\K}x\le_{\K}y$ \cite{T63}.  Hilbert metric is a distance between \emph{rays} in each part of $\K$:  $d_{H}(x,y) = 0$, $x,\,y\in\K$, $x\sim_{\K}y$, if and only if $x = \lambda y$ with $\lambda>0$.  

Hilbert and Thompson metric have been of great interest to analysts, especially for their contractivity properties. As a matter of fact, many naturally occurring maps in analysis, both linear and non-linear, are either non-expansive or contractive with respect to these metrics \cite{B73,LW94,LN12}. Moreover, it has been proven that among all projective distances $d$ on $\K$ for which the positive linear transformations are contractive w.r.t. $d$, Hilbert metric is the one with the best possible contraction ratio \cite{KP82}. 

  Thompson and Hilbert metric endow the cone with a structure of  Finsler manifold \cite{N94}. In the finite-dimensional case, the interior of the cone $\K$ defines an $n$-dimensional manifold and the tangent space at each point may be identified with $\R^{n}$. Defining the norm
\begin{equation}
\label{norm}
	\|v\|_{x}^{T} :=\inf\{\alpha>0\,:\,-\alpha x\le_{\K} v\le_{\K}\alpha x\}
\end{equation}
 on the tangent space at each point $x\in\mathring{\K}$, the length of any differentiable curve $\gamma\colon [a,b]\to\mathring{\K}$ is defined as
\[
	\ell(\gamma):=\int_{a}^{b}\|\gamma'(t)\|_{\gamma(t)}^{T}\de t.
\]
Thompson distance between any two points is recovered by minimizing over all paths connecting the points, namely
\[
	d_{T}(x,y)=\inf\{\ell(\gamma)\,:\, \gamma\in C^{1}[x,y]\},
\]
where $C^{1}[x,y]$ denotes the set of all differentiable paths $\gamma\colon [a,b]\to\mathring{\K}$ such that $\gamma(0)=x$ and $\gamma(1)=y$.  Hilbert metric is obtained along the same lines when the norm above is replaced by the semi-norm
\[
	\|v\|_x^{H} := M(v,x) - m(v,x).
\]

The Finslerian nature of Hilbert and Thompson geometries allows for the definition of minimal geodesics connecting two points in the interior of the cone $\K$. Differently from the Riemannian framework,  minimal geodesics connecting two points are usually not unique \cite[Ch. 2]{LN12}. An explicit class of minimal geodesics for the Thompson metric connecting $x,y\in\mathring{\K}$, is given by, $\chi\colon[0,1]\to\mathring{\K}$,
\begin{align}\label{eq:geod}
\chi(\tau)=\begin{cases}
\left(\frac{\beta^{\tau}-\alpha^{\tau}}{\beta-\alpha}\right)y+\left(\frac{\beta\alpha^{\tau}-\alpha\beta^{\tau}}{\beta-\alpha}\right)x, & \text{if } \beta\neq \alpha,\\
\alpha^{\tau}x, & \text{if } \beta= \alpha,
\end{cases}
\end{align}
where $\beta:=M(y,x)$ and $\alpha:=m(y,x)$. This geodesic path defines a ``projective'' straight line in the cone \cite{NW04}.

Such geodesic paths are not unique.  For instance, a distinct type of geodesic paths   connecting two positive definite matrices $X$ and $Y$ is given by 
\begin{align}\label{eq:geodT}
	\varphi_{T}(\tau)=X^{1/2}(X^{-1/2}Y X^{-1/2})^{\tau}X^{1/2}, \quad \tau\in[0,1].	
\end{align}
This path is in fact the (unique, up to a re-parametrization) Riemannian geodesic of $\S_{+}^{n\times n}$ connecting $X$ to $Y$ with respect to the affine invariant metric, see e.g. \cite[Thm. 6.1.6]{B09}.    The corresponding geodesic path in \eqref{eq:geodT}  w.r.t. Hilbert metric reads \cite[Prop. 2.6.8]{LN12}
\begin{align}\label{eq:geodH}
	\varphi_{H}(\tau)=\frac{\varphi_{T}(\tau)}{\tr(\varphi_{T}(\tau))}, \quad \tau\in[0,1],
\end{align}
where the latter path connects two unit-trace elements $X,Y\in\mathring{\S}_{+}^{n\times n}$ which are the representatives of the corresponding projective rays $\mu X$, $\mu Y$, $\mu>0$, respectively.

Finally, we remark that the Finslerian framework so far discussed for the case of finite-dimensional spaces applies without any substantial change to the case of infinite-dimensional manifolds of bounded positive self-adjoint operators on an Hilbert space. For further details on this extension we refer to the works by Corach and co-workers \cite{CPR93,CPR94}, and, in particular, to \cite{CM99,corach2000differential}.

\section{Finslerian distances in $\S_{+,\mathrm{rat}}^{n\times n}(\T)$}\label{eq:Finslerian-dist}

Since spectral densities can be thought of as bounded positive self-adjoint operators on the Hilbert space $\L_{2}^{n}[-\pi,\pi]$, the framework outlined in the previous section provides Finslerian distances in the cone $\S_{+}^{n\times n}(\T)$, and, therefore, in the space of \emph{rational} spectral densities $\S_{+,\mathrm{rat}}^{n\times n}(\T)$. Interestingly, it turns out that in the latter case the expressions of Thompson and Hilbert metric are connected with a classical problem in systems theory, the spectral factorization problem.

\begin{theorem}\label{prop1}
Consider two full normal rank spectral densities $\Phi_{1},\Phi_{2}\in \S_{+,\mathrm{rat}}^{n\times n}(\T)$ and let $W_{1},W_{2}\in\R^{n\times n}(z)$ denote the corresponding minimum-phase spectral factors. If $W_{2}^{-1}W_{1}$ has no zero/pole on $\T$, then the Hilbert and Thompson metrics between $\Phi_{1}$ and $\Phi_{2}$ are given, respectively, by
\begin{align*}
	d_{H}(\Phi_{1},\Phi_{2})&=\log\ \left\|W_{2}^{-1}W_{1}\right\|_{\H_\infty}^{2}\left\|W_{1}^{-1}W_{2}\right\|_{\H_\infty}^{2},\\ 
	d_{T}(\Phi_{1},\Phi_{2})&=\log\, \max\left\{\left\|W_{2}^{-1}W_{1}\right\|_{\H_\infty}^{2},\left\|W_{1}^{-1}W_{2}\right\|_{\H_\infty}^{2}\right\}.
\end{align*}
Otherwise, it holds $d_{H}(\Phi_{1},\Phi_{2})=d_{T}(\Phi_{1},\Phi_{2})=\infty$.
\end{theorem}
\begin{IEEEproof}
In view of the definition of $M(\cdot,\cdot)$ in \eqref{eq:M}, for any full normal rank $\Phi_{1},\Phi_{2}\in\S_{+,\mathrm{rat}}^{n\times n}(\T)$, it holds
\begin{align}
	M(\Phi_{1},\Phi_{2})&=\inf\{\lambda\in\R\,:\, \Phi_{1}(e^{j\theta})\le \lambda \Phi_{2}(e^{j\theta}), \theta\in[-\pi,\pi] \} \notag\\
	&=\inf\{\lambda\in\R\,:\, \Phi_{2}^{-\frac{1}{2}}(e^{j\theta})\Phi_{1}(e^{j\theta})\Phi_{2}^{-\frac{1}{2}}(e^{j\theta}) \le \lambda I_{n},  \theta\in[-\pi,\pi] \}\notag\\
	& =  \left\|\Phi_{2}^{-\frac{1}{2}}\Phi_{1}\Phi_{2}^{-\frac{1}{2}}\right\|_{\L_\infty},\label{eq:MPhi}
\end{align}
if $\Phi_{2}^{-\frac{1}{2}}\Phi_{1}\Phi_{2}^{-\frac{1}{2}}$ is analytic on $\T$, and $M(\Phi_{1},\Phi_{2})=\infty$ otherwise.
In order to deal with rational matrix-valued functions we can replace, without affecting the value of $M(\Phi_{1},\Phi_{2})$, the square root $\Phi_{2}^{1/2}$ in the latter expression with the minimum-phase spectral factor $W_{2}\in\R^{n\times n}(z)$ of $\Phi_{2}$. (In fact, $\Phi_{2}^{\frac{1}{2}}U=W_{2}$, where $U$ is a suitable $n\times n$ unitary matrix-valued function on $\T$). Therefore, Equation~\eqref{eq:MPhi} becomes
\begin{align*}
	M(\Phi_{1},\Phi_{2}) &= \left\|\Phi_{2}^{-\frac{1}{2}}\Phi_{1}\Phi_{2}^{-\frac{1}{2}}\right\|_{\L_\infty}\\
	&= \left\|W_{2}^{-1}\Phi_{1}W_{2}^{-*}\right\|_{\L_\infty}\\
	&=\left\|W_{2}^{-1}W_{1}W_{1}^{*}W_{2}^{-*}\right\|_{\L_\infty}\\
	&=\left\|W_{2}^{-1}W_{1}\right\|^{2}_{\L_\infty},
\end{align*}
if $W_{2}^{-1}W_{1}$ has no pole on $\T$, and $M(\Phi_{1},\Phi_{2})=\infty$ otherwise.
Further, if $W_{2}^{-1}W_{1}$ has no pole on $\T$, $W_{2}^{-1}W_{1}$ is analytic in (an open set containing) the complement of $\D$, so that we have
\[
	M(\Phi_{1},\Phi_{2})=\left\|W_{2}^{-1}W_{1}\right\|^{2}_{\H_{\infty}},
\]
where we have replaced the $\L_\infty$-norm with the $\H_{\infty}$-norm. Similarly, we have 
\[
	M(\Phi_{2},\Phi_{1})=\left\|W_{1}^{-1}W_{2}\right\|^{2}_{\H_{\infty}},
\]
if $W_{1}^{-1}W_{2}$ has no pole on $\T$, or, equivalently, if $W_{2}^{-1}W_{1}$ has no zero on $\T$, and $M(\Phi_{2},\Phi_{1})=\infty$ otherwise.

Eventually, observing that $m(\Phi_{1},\Phi_{2})=M(\Phi_{2},\Phi_{1})^{-1}$, a substitution of the values of $M(\Phi_{1},\Phi_{2})$ and $M(\Phi_{2},\Phi_{1})$ into the expressions of Hilbert and Thompson metrics in Definition \ref{def1} yields the thesis.
\end{IEEEproof}

\begin{remark}
The proof of Theorem \ref{prop1} shows that  the expressions of the Hilbert and Thompson metric still hold if we replace the canonical (i.e., minimum-phase) spectral factors of the two spectra $\Phi_{1},\,\Phi_{2}$ with any other spectral factor of $\Phi_{1},\,\Phi_{2}$ (i.e., spectral factors not necessarily analytic in the complement of $\D$ and with analytic inverse in the complement of the closure of $\D$). The important difference is that, in this case, the $\H_\infty$-norm must be replaced by the $\L_\infty$-norm.
\end{remark}

\begin{remark}
As discussed in the previous section, the difference between Hilbert and Thompson metric consists of the fact that the Thompson metric is a bona fide distance on each part of $\S_{+,\mathrm{rat}}^{n\times n}(\T)$ (and, in particular, on its interior), while Hilbert metric is a distance between rays in each part of the latter cone. It is worth remarking that projective invariance has proved to be a desirable property since in many applications, such as spectral estimation or speech processing, the shape of the spectral densities rather than their relative scalings is the discriminative feature \cite{AV14,G07}.
\end{remark}
 
 \begin{remark}\label{rem:art}
The expressions of Hilbert and Thompson metrics in Theorem \ref{prop1} apply also to the case of general non-rational spectral densities in $\S_{+}^{n\times n}(\T)$.\footnote{Notice that in case the minimum-phase spectral factors of $\Phi_{1}$, $\Phi_{2}$ do not exist, the expressions in Theorem \ref{prop1} still holds by replacing the minimum-phase spectral factors with the corresponding frequency-wise matrix square roots $\Phi_{1}^{1/2}$, $\Phi_{2}^{1/2}$.}
In this case, however,  one issue that arises is that the distance between almost identical  spectral densities can be made arbitrarily large. With reference to the scalar case, this occurs when one of the two spectral densities exhibits a very sharp and narrow frequency peak. For the sake of illustration, consider the two scalar spectral densities in $\S_{+}^{1\times 1}(\T)$
\begin{align}\label{eq:remark_path}
	\phi_1(e^{j\theta}) = 1, \quad \phi_{2,\varepsilon}(e^{j\theta}) = \begin{cases} \varepsilon^{-1} & \text{ if } |\theta|\,\leq\varepsilon, \\ 1 & \text{ otherwise,} \end{cases}
\end{align}
where $\theta\in[-\pi,\pi]$ and $\varepsilon>0$. It can be seen that, for $\varepsilon\to 0$, $d_H(\phi_1,\phi_{2,\varepsilon})\to\infty$ and $d_T(\phi_1,\phi_{2,\varepsilon})\to\infty$, in spite of the fact that the two spectral densities are identical with the only exception of a neighborhood of the frequency $\theta=0$ (see also Figure \ref{Fig:remark}). Importantly, when restricting the attention to spectral densities that are  ``sufficiently regular'', e.g., those belonging to the space of rational spectral densities with bounded McMillan degree, these pathological cases are ruled out.

In view of the above property, the proposed distances are not suitable to treat spectral densities featuring ``spectral lines''.  On the other hand,  the same property is a relevant feature 
in a distance  for the THREE spectral estimation problem described in Example \ref{ex:THREE-est}. Indeed, in view of this property the THREE problem with a distance as the ones introduced in this paper cannot produce solutions featuring ``artifacts'' (spurious spectral lines) which is one of the main drawbacks of the method with classical distances. In this sense the distances just introduced can be used for a   ``robust'' version of the THREE problem (see also the concluding remarks for more details).

\begin{figure}[!h]
\begin{center}
\begin{tikzpicture}
\draw[-latex] (-3,0) -- (3,0) node[right] {\footnotesize $\theta$};
\draw[-latex] (0,-0.5) -- (0,3);
\draw[dotted] (-2.75,-0.2) node[below] {\footnotesize $-\pi$} -- (-2.75,2.75);
\draw[dotted] (2.75,-0.2) node[below=0.05cm] {\footnotesize $\pi$} -- (2.75,2.75);
\draw[red,thick] (-2.75,0.5) -- (2.75,0.5);
\draw[blue,thick,dashed] (-2.75,0.5) -- (-0.75,0.5);
\draw[blue,thick,dashed] (0.75,0.5) -- (2.75,0.5);
\draw[blue,thick,dashed] (0.75,0.5) -- (2.75,0.5);
\draw[blue,thick,dashed] (-0.75,1.5) -- (0.75,1.5);
\draw[blue,dotted] (-0.75,1.5) -- (-0.75,0.5);
\draw[blue,dotted] (0.75,1.5) -- (0.75,0.5);
\draw[ForestGreen,thick,densely dotted] (-2.75,0.5) -- (-0.25,0.5);
\draw[ForestGreen,thick,densely dotted] (0.25,0.5) -- (2.75,0.5);
\draw[ForestGreen,thick,densely dotted] (0.25,0.5) -- (2.75,0.5);
\draw[ForestGreen,thick,densely dotted] (-0.25,2.5) -- (0.25,2.5);
\draw[ForestGreen,dotted] (-0.25,2.5) -- (-0.25,0.5);
\draw[ForestGreen,dotted] (0.25,2.5) -- (0.25,0.5);
\draw[dotted] (-0.25,0.5) -- (-0.25,-0.25) node[below] {\footnotesize $-\varepsilon_1$};
\draw[dotted] (0.25,0.5) -- (0.25,-0.25) node[below=0.05cm] {\footnotesize $\varepsilon_1$};
\draw[dotted] (-0.75,0.5) -- (-0.75,-0.25) node[below] {\footnotesize $-\varepsilon_2$};
\draw[dotted] (0.75,0.5) -- (0.75,-0.25) node[below=0.05cm] {\footnotesize $\varepsilon_2$};
\node at (0.125,0.7) {\footnotesize $1$};
\node at (0.325,2.75) {\footnotesize $\varepsilon_1^{-1}$};
\node at (0.325,1.75) {\footnotesize $\varepsilon_2^{-1}$};
\draw[red,thick] (1.3,2.5) -- (1.8,2.5) node[right] {\footnotesize $\phi_1$};
\draw[blue,thick,dashed] (1.3,2.15) -- (1.8,2.15) node[right] {\footnotesize $\phi_{2,\varepsilon_{2}}$};
\draw[ForestGreen,thick,densely dotted] (1.3,1.8) -- (1.8,1.8) node[right] {\footnotesize $\phi_{2,\varepsilon_{1}}$};
\end{tikzpicture}
\caption{Qualitative plot of the spectral densities in \eqref{eq:remark_path} for two values of $\varepsilon$, i.e. $\varepsilon_2>\varepsilon_1>0$.}
\label{Fig:remark}
\end{center}
\end{figure}
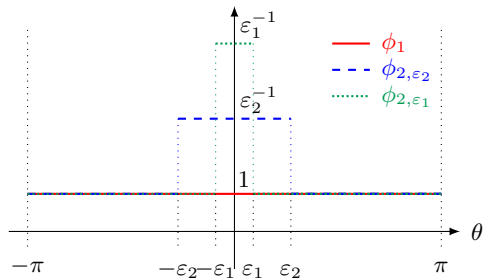
\end{remark}


\subsection{Filtering invariance}\label{subsec:invar}

The above-introduced distances possess the following important property:   
\begin{equation}
\label{filterinvariance}
\forall\, T\in\R^{n\times n}_*(z): \; \;  d (\Phi_{1},\Phi_{2})=d (T\Phi_{1}T^{*},T\Phi_{2}T^{*}).
\end{equation}

This property readily follows from the definition of  Hilbert and Thompson distances and the expression of $M(\Phi_{1},\Phi_{2})$ in \eqref{eq:MPhi}.  Since the set $\R^{n\times n}_*(z)$ defines a group, the mapping $\Phi \mapsto T\Phi T^{*}$ defines a congruence group action of $\R^{n\times n}_*(z)$ on the set of rational spectral densities. This group action is transitive, that is, any rational spectral density can be obtained by acting on the identity element. 

A metric that satisfies (\ref{filterinvariance}) can be said to be {\it filtering invariant} because of the following statistical interpretation. Any spectral density $\Phi$ with minimum-phase spectral factor $W$ can be identified to a $n$-dimensional zero-mean second-order stationary purely nondeterministic stochastic process $\{y(t)\}_{t\in\Z}$
generated by filtering a white noise process through $W$.  The action $\Phi \mapsto T\Phi T^{*}$ has therefore the interpretation of filtering the process
 with the linear time-invariant filter $T\in\R_{*}^{n\times n}(z)$. Likewise, the property (\ref{filterinvariance}) has the interpretation that the distance between two spectral densities, or, equivalently, two zero-mean second-order stationary purely nondeterministic stochastic processes, is unchanged when the two processes are filtered by the same filter.

Any filtering invariant metric  is entirely specified by defining the distance to identity. Furthermore, one has $d(\Phi,I)=d(\Phi^{-1},I)$. In other words, the distance is a {\it distortion} measure.  


Filtering invariance is a fundamental property of  classical metrics. In the scalar case, $\phi_{1},\phi_{2}\in\S_{+,\mathrm{rat}}^{1\times 1}(\T)$, the log spectral deviation \cite{GBGM80}
$$
\left(\int_{-\pi}^\pi\left\| \log \frac{\phi_1}{\phi_2}  \right\|^2\frac{\de \theta}{2\pi}\right)^{1/2} $$
 is an early example of filtering invariant distortion measure.
The recent work \cite{JNG12} shows that the multivariate generalization  
\begin{align}\label{eq:Riemann}
	d_R(\Phi_1,\Phi_2)&=\left(\int_{-\pi}^\pi\left\| \log \Phi_1^{-1/2}\Phi_2\Phi_1^{-1/2} \right\|_F^2\frac{\de \theta}{2\pi}\right)^{1/2}\notag\\
	&=\left(\int_{-\pi}^\pi\left\| \log W_1^{-1}\Phi_2W_1^{-*} \right\|_F^2\frac{\de \theta}{2\pi}\right)^{1/2}
\end{align}
is the unique Riemannian bona fide distance that is filtering invariant. This metric is a natural generalization of the affine invariant metric between positive definite matrices. Affine invariance corresponds to filtering invariance in the static case: the congruence group action  reduces to an action of the general linear group. The metric is in this case a distance between  $n$-dimensional zero-mean second-order random vectors, and the invariance property is an invariance with respect to an affine change of coordinates. The importance of this invariance property in the context of estimation problems has been emphasized for instance in \cite{smith2005}. In \cite{JNG12}, filtering invariance emerges as a natural property when measuring the ``flatness'' of innovations processes.
Filtering invariance is also a leading prerequisite in the work of Martin \cite{martin2000}, whose resulting metric, which applies to scalar spectral densities $\phi_{1},\phi_{2}\in\S_{+,\mathrm{rat}}^{1\times 1}(\T)$, can be written as 
\begin{align}\label{eq:Martin}
	d_M(\phi_1,\phi_2)&=\left(\int_{-\pi}^\pi\left(\mathfrak{D}^{\frac{1}{2}}\log\frac{\phi_{1}}{\phi_{2}}\right)^{2}\frac{\de \theta}{2\pi}\right)^{1/2},
\end{align} 
where $\mathfrak{D}^{\lambda}$, $\lambda>0$, is the fractional derivative operator in the frequency domain. Finally, filtering invariance is also a key property of
the classical Itakura--Saito divergence (see \cite{IS70,FMP12}).

The  Riemannian distance (\ref{eq:Riemann}) and the Thompson metric introduced in Section 2  are thus close relatives: they are  bona fide distances which satisfy filtering invariance and endow the cone of spectral densities with a differential metric structure. The first one induces a Riemannian structure through  an invariant inner product (that reduces to the standard inner product at identity),  while the second induces a Finslerian structure through the invariant norm (\ref{norm}). Both distances depend on the same log spectral quantity 
frequency-wise, but the Riemannian distance results in a two-norm of that frequency-domain function, whereas the Finsler distance results in an infinite-norm.

\subsection{Computational properties}

Theorem \ref{prop1} shows that the computation of Hilbert and Thompson metrics in the cone of rational spectral densities essentially consists of: (i) the calculation of the minimum-phase spectral factors $W_{1}$ and $W_{2}$, and (ii) the calculation of the $\H_{\infty}$-norm of the ``ratio'' of the latter spectral factors. 
For the sake of clarity, we illustrate below a simple example of computation of these metrics.
 \begin{example}
 Consider the scalar rational spectral densities
 \begin{align*}
 	\phi_{1}(z) &= -\frac{4}{2z+-5+2z^{-1}},  \\
	\phi_{2}(z) &= \frac{9}{3z+10+3z^{-1}}.
 \end{align*}
 The minimum-phase spectral factors of $\phi_{1}$ and $\phi_{2}$ are given by, respectively,
  \begin{align*}
 	w_{1}(z) = \frac{z}{z-1/2}, \ \ \ w_{2}(z) = \frac{z}{z+1/3}.
 \end{align*}
 By direct computation, we have 
   \begin{align*}
 	\left\|\frac{w_{2}}{w_{1}}\right\|_{\mathcal{H}_{\infty}}^{2} &=\sup_{|z|> 1}\frac{|z+1/3|}{|z-1/2|}= \frac{8}{3}  \\
	\left\|\frac{w_{1}}{w_{2}}\right\|_{\mathcal{H}_{\infty}}^{2} &=\sup_{|z|> 1}\frac{|z-1/2|}{|z+1/3|}= \frac{9}{4}. 
 \end{align*}
 Hence, in view of Theorem \ref{prop1},
 \begin{align*}
	d_{H}(\phi_{1},\phi_{2})&=\log\ \left\|\frac{w_{2}}{w_{1}}\right\|_{\mathcal{H}_{\infty}}^{2}\left\|\frac{w_{1}}{w_{2}}\right\|_{\mathcal{H}_{\infty}}^{2}=\log 6,\\ 
	d_{T}(\phi_{1},\phi_{2})&=\log\, \max\left\{\left\|\frac{w_{2}}{w_{1}}\right\|_{\mathcal{H}_{\infty}}^{2},\left\|\frac{w_{1}}{w_{2}}\right\|_{\mathcal{H}_{\infty}}^{2}\right\} =\log \frac{8}{3}.
\end{align*}
 \end{example} \vspace{0.15cm}

Remarkably, the computation of (minimum-phase) spectral factors and $\H_{\infty}$-norms represent two extensively studied problems in systems and control theory and several algorithms are available in the literature to perform these steps, even in the most challenging multivariate case. More specifically, a general result on the existence and (essential) uniqueness of minimum-phase spectral factors can be found  
in \cite{youla1961factorization} for the continuous-time case and in \cite{baggio2016factorization} for the discrete-time case. An algorithm to compute this spectral factor is described there, however, for the calculation of minimum-phase spectral factors there exist several more efficient routines based on the solution of suitable Stein and Riccati Equations, see e.g. \cite{O05}. 
 Whereas, an efficient method for computing $\H_\infty/\L_\infty$-norms of rational matrix-valued transfer functions is the Boyd--Balakrishnan--Bruinsma--Steinbuch method \cite{BS90,BK90} which is based on an iterative bisection-like algorithm and leads to quadratic convergence.\footnote{The approach is formulated in the continuous-time case, however there exist variants of this algorithm for computing the discrete-time $\H_{\infty}$-norm, based on computing eigenvalues of symplectic instead of Hamiltonian matrices, see e.g. \cite{GVV98}.} Moreover, an upper bound to these norms can be found by inspecting the eigenvalues of the symplectic matrix associated with the state-space representation of the system \cite[Lemma 21.10]{ZDG96}.

In view of the previous analysis, a merit of the Finslerian (Thompson) distance over its Riemannian relative \eqref{eq:Riemann} is at the  computational level. In fact, the calculation of the Riemannian distance requires the frequency-wise computation of the (matrix) logarithm of $W_1^{-1}\Phi_2W_1^{-*}$, an operation which appears numerically challenging in the multivariate setting.  In contrast, the calculation of the Thompson metric involves the computation of minimum-phase spectral factors and $\H_{\infty}$-norms, for which  efficient numerical algorithms are available.

For completeness, it should be mentioned that one way of overcoming the computational burden of the Riemannian distance is to replace it with a divergence measure. In the static case, Kullback--Leibler divergence  approximates the Riemannian distance up to third order. In the dynamic case, the paper \cite{JNG12} considers quadratic approximations of divergence measures. In the rational case, one such quantity takes the form
\begin{align*}
	d_{F}(\Phi_{1},\Phi_{2})=\|W_2^{-1} W_1\|_{\H_2}^{2}+\|W_1^{-1} W_2\|_{\H_2}^{2}-2n,
\end{align*} 
where $\|\cdot\|_{\H_{2}}$ denotes the $\H_{2}$-norm of a discrete-time transfer function \cite[Sec.~4.3]{ZDG96}.
It is not a distance (in fact, it does not obey the triangle inequality) but it provides a tractable quadratic approximation of the Riemannian distance.

\subsection{Geodesic paths}

One can generalize the geodesic expression \eqref{eq:geodT} obtained for the positive definite matrix case to the case of bounded positive operators on an Hilbert space \cite{CPR93,CPR94}. From this generalization, it follows that a minimal geodesic path in $\S_{+}^{n\times n}(\T)$ connecting $\Phi_{1},\,\Phi_{2}\in \mathring{\S}_{+,\mathrm{rat}}^{n\times n}(\T)$ w.r.t. Thompson metric is given by\footnote{Notice that we have expressed the geodesic in terms of the rational spectral factor $\smash{W_1}$ instead of the frequency-wise matrix square root $\smash{\Phi_{1}^{1/2}}$. This equivalent rewriting follows from the fact that $\smash{\Phi_{1}^{1/2}=W_{1}U}$ with $U$ being an $n\times n$ unitary matrix-valued function on $\T$.}
\begin{align}\label{eq:gT}
\varphi_{T}(\tau)=W_1(W_1^{-1}\Phi_{2} W_1^{-*})^{\tau}W_1^*, \quad \tau\in[0,1],
\end{align}
where $W_{1}\in\R^{n\times n}(z)$ is the minimum-phase spectral factor of $\Phi_{1}$. Notice that Equation \eqref{eq:gT} coincides with the (unique, up to a re-parametrization) Riemannian geodesic between spectral densities \cite{JNG12}.

However, as discussed in Sec.~\ref{sec:dist-cones}, the Finslerian framework allows for the definition of multiple minimal geodesic paths. For instance, an alternative explicit Thompson's minimal geodesic for spectral densities follows the ``projective'' straight line interpolation in \eqref{eq:geod}. This yields the following minimal geodesic path connecting $\Phi_{1},\,\Phi_{2}\in \mathring{\S}_{+,\mathrm{rat}}^{n\times n}(\T)$ w.r.t. Thompson metric
\begin{align}\label{eq:cT}
\chi_{T}(\tau)=\begin{cases}
\left(\frac{\beta^{\tau}-\alpha^{\tau}}{\beta-\alpha}\right)\Phi_{2}+\left(\frac{\beta\alpha^{\tau}-\alpha\beta^{\tau}}{\beta-\alpha}\right)\Phi_{1}, & \text{if } \beta\neq \alpha,\\
\alpha^{\tau}\Phi_{1}, & \text{if } \beta= \alpha,
\end{cases}
\end{align}
where $\tau\in[0,1]$, 
\[
\beta:=\left\|W_{1}^{-1}W_{2}\right\|^{2}_{\H_{\infty}}, \quad \alpha:=1/\left\|W_{2}^{-1}W_{1}\right\|^{2}_{\H_{\infty}},\] 
with $W_{1}\in\R^{n\times n}(z)$ and $W_{1}\in\R^{n\times n}(z)$ being the minimum-phase spectral factors of $\Phi_{1}$ and $\Phi_{2}$, respectively. 

When applied to the whole cone $\mathring{\S}_{+}^{n\times n}(\T)$, the Riemannian geodesic path \eqref{eq:gT} renders this space geodesically complete, meaning that for all $\tau\in\R$, $\varphi_{T}(\tau)$ belongs to  $\mathring{\S}_{+}^{n\times n}(\T)$. This is indeed a remarkable property that allows for extrapolating along the geodesic paths. 
On the other hand, when considering rational spectral densities, $\varphi_{T}(\tau)$ is not, in general, rational. Concerning the Finslerian geodesic path in \eqref{eq:cT} we have the following result.

\begin{proposition}\label{prop3}
For all $\Phi_{1},\Phi_{2}\in\mathring{\S}_{+,\mathrm{rat}}^{n\times n}(\T)$ and $\tau\in\R$, $\chi_{T}(\tau) \in\mathring{\S}_{+,\mathrm{rat}}^{n\times n}(\T)$.
\end{proposition}
\begin{IEEEproof} 
The path $\chi_{T}(\tau)$ is bounded and rational for all $\tau\in\R$ by construction since $\alpha$ and $\beta$ are bounded real scalars. So it remains to prove positivity of $\chi_{T}(\tau)$. The case $\beta =\alpha$ is straightforward, so in what follows we suppose $\beta \ne \alpha$. Notice that, by pre- and post-multiplying $\chi_{T}(\tau)$ by $\Phi^{-1/2}$ and then diagonalizing the resulting expression, the condition $\chi_{T}(\tau)>0$ for all $\theta\in[-\pi,\pi]$ can be seen to be equivalent to
	\begin{align}\label{eq:proof-geod}
\left(\frac{\beta^{\tau}-\alpha^{\tau}}{\beta-\alpha}\right)\Lambda(e^{j\theta})+\left(\frac{\beta\alpha^{\tau}-\alpha\beta^{\tau}}{\beta-\alpha}\right)I>0,
\end{align}
for all $\theta\in[-\pi,\pi]$, where 
\[
	\Lambda(e^{j\theta}):=\diag[\lambda_{1}(e^{j\theta}),\lambda_{2}(e^{j\theta}),\dots,\lambda_{n}(e^{j\theta})],
\]
$\lambda_{1}(e^{j\theta})\ge\lambda_{2}(e^{j\theta})\ge \dots\ge \lambda_{n}(e^{j\theta})>0$, $\theta\in[-\pi,\pi]$, has in its diagonal the frequency-wise eigenvalues of $\Phi_{1}^{-1/2}(e^{j\theta})\Phi_{2}(e^{j\theta})\Phi_{1}^{-1/2}(e^{j\theta})$. Next, we note that 
\[
	\lambda_{n}(e^{j\theta})\ge \min_{\theta\in[-\pi,\pi]} \lambda_{n}(e^{j\theta}) = 1/\left\|W_{2}^{-1}W_{1}\right\|^{2}_{\H_{\infty}}=\alpha.
\]
In view of Eq.~\eqref{eq:proof-geod}, this in turn implies that for all $\theta\in[-\pi,\pi]$,
\begin{align*}
&\left(\frac{\beta^{\tau}-\alpha^{\tau}}{\beta-\alpha}\right)\Lambda(e^{j\theta})+\left(\frac{\beta\alpha^{\tau}-\alpha\beta^{\tau}}{\beta-\alpha}\right)I\\
&\ge\left(\frac{\beta^{\tau}-\alpha^{\tau}}{\beta-\alpha}\right)\alpha I+\left(\frac{\beta\alpha^{\tau}-\alpha\beta^{\tau}}{\beta-\alpha}\right)I =\alpha^{t}I>0,
\end{align*}
which completes the proof.
\end{IEEEproof}

In more practical terms, the previous proposition states that one can  interpolate along the geodesic path \eqref{eq:cT} while remaining in the cone of \emph{rational} spectral densities. This is a fundamental feature of Finsler geometry that does not have a Riemannian counterpart. In light of this fact, we argue that, when dealing with rational spectral densities, the Finsler geodesic \eqref{eq:cT} may be a more natural choice when compared to the Riemannian geodesic \eqref{eq:gT}.
In view of the geodesical completeness, this is true also for {\em extrapolation} so that,
given two rational spectral densities  $\Phi_{1}$ and $\Phi_{2}$,  
we can select a rational spectral density in the geodesic line connecting $\Phi_{1}$ and $\Phi_{2}$ and this spectral density is not necessarily between $\Phi_{1}$ and $\Phi_{2}$ but may also be  chosen to be ``before'' $\Phi_{1}$ or ``after'' $\Phi_{2}$. This is particularly interesting in applications as illustrated in the next section where the morphing between a male and a female voice is discussed and we can, for example, go ``beyond male'' and synthesize a particularly baritonal voice.


\begin{remark}
Geodesic expressions \eqref{eq:gT} and \eqref{eq:cT} for Thompson metric can be adapted to geodesics for Hilbert metric by considering the corresponding ``normalized'' versions. For instance, w.r.t. ``normalized'' spectra $\Phi_{1}$, $\Phi_{2}\in \mathring{\S}_{+,\mathrm{rat}}^{n\times n}(\T)$ such that $\int_{-\pi}^{\pi} \tr(\Phi_{1})\frac{\de\theta}{2\pi}=\int_{-\pi}^{\pi} \tr(\Phi_{2})\frac{\de\theta}{2\pi}=1$, geodesic \eqref{eq:gT} becomes 
\begin{align}\label{eq:gH}
\varphi_{H}(\tau)=\frac{W_1(W_1^{-1}\Phi_{2} W_1^{-*})^{\tau}W_1^*}{\int_{-\pi}^{\pi}\tr(W_1(W_1^{-1}\Phi_{2} W_1^{-*})^{\tau}W_1^*)\frac{\de\theta}{2\pi}}, \quad \tau\in[0,1].
\end{align}
\end{remark}

\begin{remark}
The reader will observe that the geodesic paths discussed in this section inherit the invariance properties of the metric discussed in Subsection~\ref{subsec:invar}. Hence, the construction of geodesic connecting curves between any two points can always be recast as the construction of a geodesic connecting curve between an arbitrary point and the identity.
\end{remark}

\section{Applications: Speech morphing}\label{sec:appl}

In this section, we show how to apply the geometric tools we developed in the previous section as a means to interpolate and extrapolate rational spectra describing the frequency content of speech data. More precisely, we analyze the task of morphing the voice of an individual into the voice of another individual. 
To this end, we first review some standard facts concerning speech modelling and synthesis, that can be found, for instance, in  \cite{RS78,MG82,HAHR01}.

Speech signals can be considered approximately stationary when restricted to a small time interval (typically, $\sim 25$ ms). Within such an interval, a speech signal can be modelled by a linear time-invariant filter driven by a suitable excitation signal. 
For each time fragment, the morphing of two speech signals can be accomplished via interpolation of the  (rational) spectral densities describing the two modelled signals. Starting from the morphed spectrum, a morphed speech signal can be generated as an output of suitable linear filter, similarly as before. The complete morphed speech signal can be eventually recovered by stacking together all the morphed audio fragments.

As suggested in \cite{JTG08} and briefly mentioned in Section~\ref{sec:motivations}, the task of interpolating two spectral densities can be naturally and efficiently carried out using suitable geodesic paths. Here, we focus on spectral interpolation via the Finslerian geodesic in \eqref{eq:cT}.

Before illustrating the obtained results, we briefly discuss the implementation details of the morphing procedure and how the latter compares to the method proposed in \cite{JTG08}. 

\begin{figure}[!h]
\begin{center}
\includegraphics[scale=0.65]{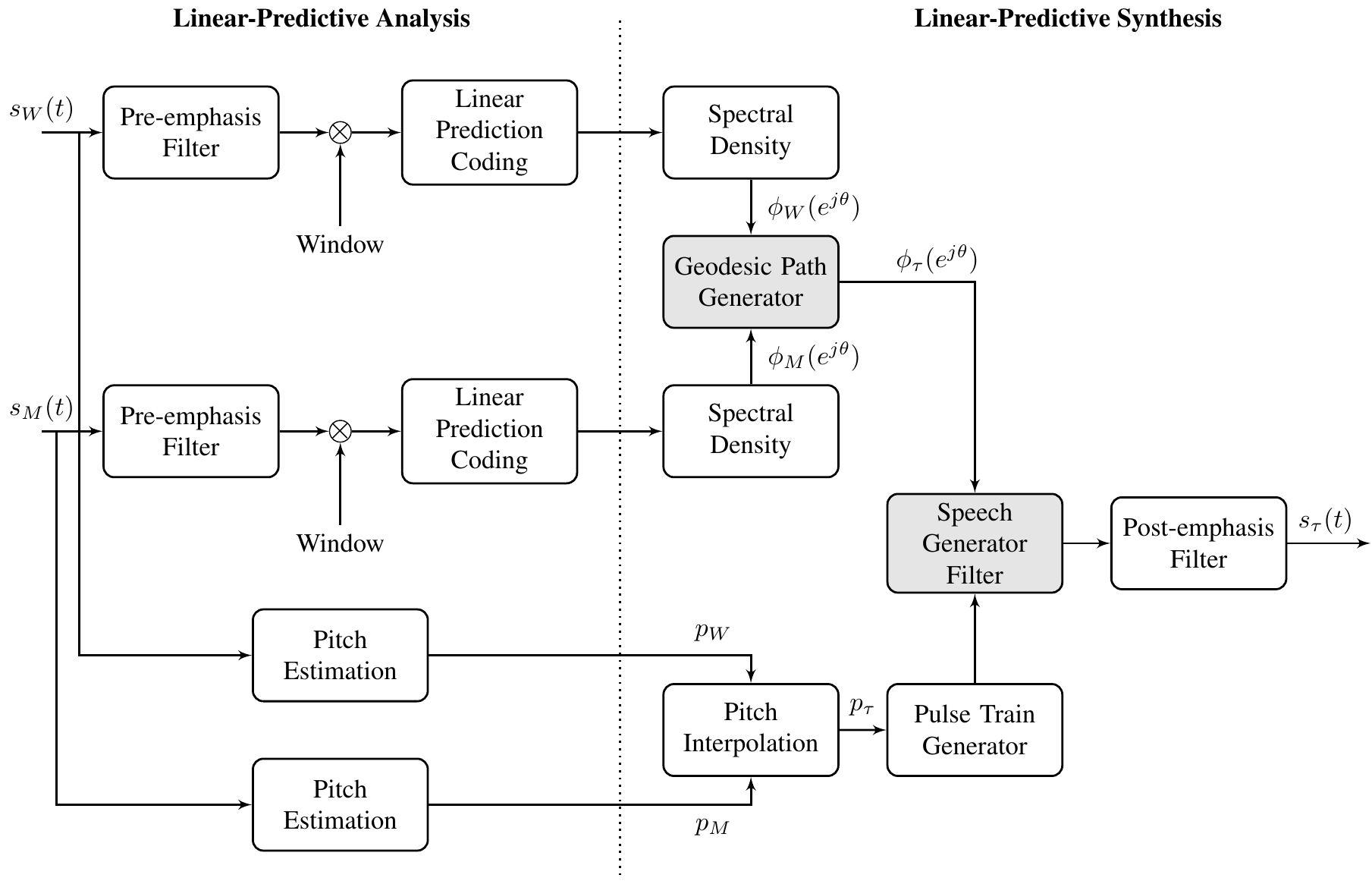}
\caption{Block diagram of the speech morphing procedure. The scheme is divided in two main blocks: Linear-Predicitive Analysis wherein a linear AR model for each speech signal is estimated, and Linear-Predicitive Synthesis wherein the morphed speech is synthesized. The main differences w.r.t.~the morphing method proposed in \cite{JTG08} are highlighted in gray.}
\label{Fig:s1}
\end{center}
\end{figure}

The morphing approach we considered is schematically depicted in the block diagram of  Figure~\ref{Fig:s1}. 
Here, we consider two audio samples $s_{M}(t)$ and $s_{W}(t)$ corresponding to the phoneme /\textscripta\textlengthmark/ spoken by a male and a female individual, respectively. Each audio signal is sampled at 16 kHZ and has a duration of $0.3$ s. We partition the signals into frames of $25$ ms, and we estimate the pitch period of the male speech signal ($p_{M}$) and of the female speech signal ($p_{W}$) via residual-based estimation \cite[Sec.~6.7]{HAHR01}. As common practice, signals $s_{M}(t)$ and $s_{W}(t)$ are pre-filtered with a ``pre-emphasis'' filter in order to reduce their low-frequency content. Then an Hamming window convolution is applied to each filtered signal. The linear model estimation (known also as linear predictive coding) is performed using classical techniques: first we estimate the covariance lags of the signals via auto-correlation method and then we apply the Levinson--Durbin method to obtain the AR model coefficients\footnote{Here, we fixed the AR model order to 14.} \cite[Sec.~6.3]{HAHR01}. From the linear model the corresponding rational power spectral densities of the estimated signals, denoted by $\phi_{M}(e^{j\theta})$ and $\phi_{W}(e^{j\theta})$, are computed and then interpolated using Finsler geodesic  \eqref{eq:cT}. The resulting rational morphed spectrum is denoted by $\phi_{\tau}(e^{j\theta})$, $\tau\in[0,1]$, where $\phi_{0}:= \phi_{W}$ and $\phi_{1}:= \phi_{M}$. The synthesis of the morphed speech signal is then simply performed by feeding the minimum-phase rational spectral factor of the morphed spectrum $\phi_{\tau}(e^{j\theta})$ with an excitation signal consisting of a pulse train with frequency $p_{\tau}$, where $p_{\tau}$ is obtained by linear or geometric interpolation of $p_{M}$ and $p_{W}$. Finally, a ``post-emphasis'' filter is applied to the resulting speech signal in order to compensate the effect of the ``pre-emphasis'' filter.

The above illustrated morphing scheme shares many similarities with the procedure adopted in \cite{JTG08}. However, a main difference is apparent. In \cite{JTG08} the synthesis step requires the estimation of a (AR) linear model from the (typically non-rational) morphed spectrum. Using Finslerian geodesic \eqref{eq:cT} this additional step can be bypassed since $\phi_{\tau}(e^{j\theta})$ is rational by construction, allowing for a (considerable) simplification of the digital implementation procedure.

\begin{figure}[h!]
\begin{center}
\includegraphics[scale=0.75]{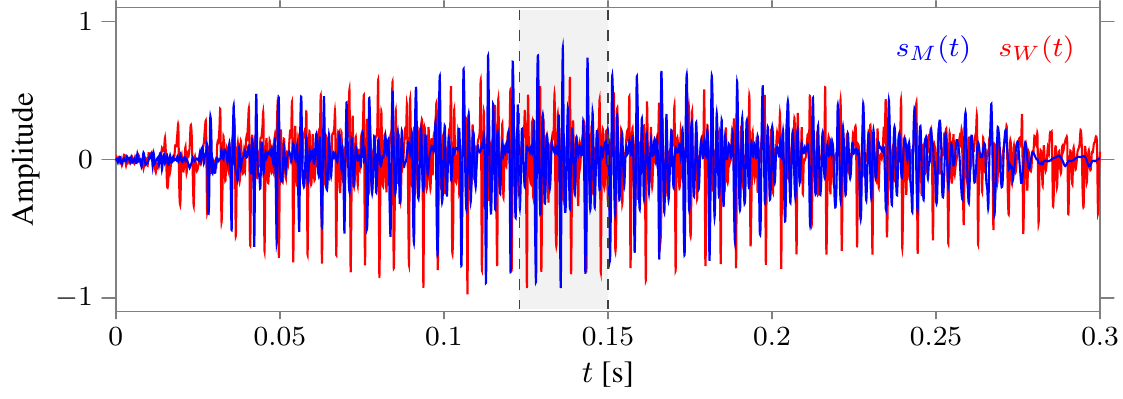}

\vspace{0.25cm}

\includegraphics[scale=0.75]{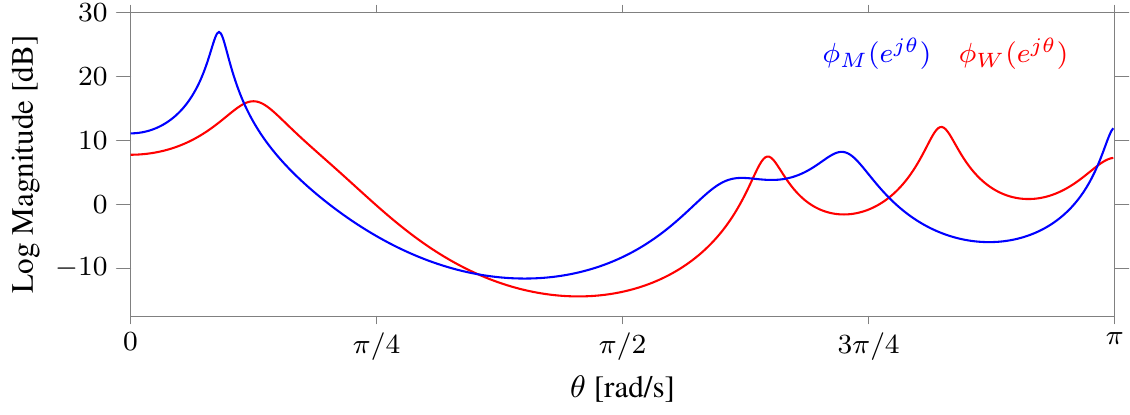}
\caption{Top plot: Speech signals corresponding to phoneme  /\textscripta\textlengthmark/ spoken by a male ($s_{M}(t)$) and female individual ($s_{W}(t)$). Bottom plot: (One-side) rational spectral densities obtained from linear predictive modelling of $s_{M}(t)$ and $s_{W}(t)$ within the time interval highlighted in the dashed gray zone in the top plot.}
\label{Fig:s2}
\end{center}
\end{figure}

The top plot of Figure~\ref{Fig:s2} shows the temporal behavior of the two considered speech waveforms $s_{M}(t)$ and $s_{W}(t)$. The dashed gray zone in the top plot highlights a specific time fragment. With reference to this time fragment, the bottom plot of Figure~\ref{Fig:s2} shows the two rational spectral densities estimated using the linear predictive estimation procedure outlined above.

\begin{figure}[!h]
\begin{center}
\includegraphics[scale=0.875]{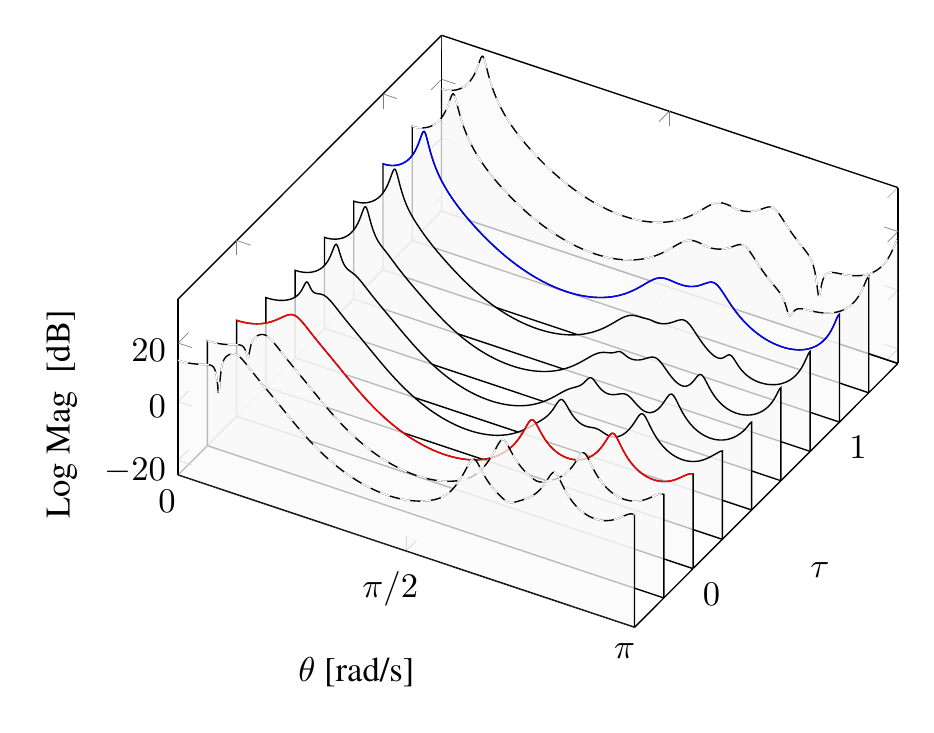}
\end{center}
\caption{Interpolation of $\phi_{M}(e^{j\theta})$ (blue curve) and $\phi_{W}(e^{j\theta})$ (red curve) via the Finslerian geodesic \eqref{eq:cT} for different values of $\tau\in\R$. The dashed curves denote extrapolated spectra. }
\label{Fig:s3}
\end{figure}

In Figure~\ref{Fig:s3}, the result of interpolation between $\phi_{M}(e^{j\theta})$ and $\phi_{W}(e^{j\theta})$ via Finslerian geodesic \eqref{eq:cT} is depicted. The dashed curves denote spectral densities obtained by extrapolating along the geodesic path. From this figure, it is interesting to observe that the interpolation/extrapolation behavior of the peaks of the morphed spectral density seems to be almost linear in a logarithmic scale. 

Results of morphing for full words and sentences are available in an audible format in \cite{speech}. In spite of the quite rudimentary nature of the procedure and the fact that the effect of the acoustic quality is very subjective, the obtained morphed audio signal seems often to be of acceptable quality. Of course, the proposed morphed procedure can be thought of as a building block which allows for the incorporation of more sophisticated processing tools to enhance the quality of the resulting speech signal.

\section{Concluding remarks and future directions}\label{sec:conclusions}

In this paper, we studied a class of conal distances for rational spectral densities arising from Finsler geometry. The proposed Thompson and Hilbert metric have a number of attractive and unique features. In particular, they are easily computable, they enjoy filtering invariance and they possess minimal geodesics that preserve rationality. These properties make these distances suitable for application in a variety of problems across systems and control theory and signal processing.

For instance, a problem that could benefit from the use of such distances is spectral estimation in the THREE-like framework discussed in the introductory Example~\ref{ex:THREE-est}. A common feature of all  distances proposed so far to tackle this problem is that they involve the two-norm of a frequency-wise quantity defined on the unit circle.  Choosing the Finsler distances of this paper could lead to a  \emph{robust} version of THREE-like spectral estimation.   A main motivation for this modified formulation concerns the reduction of \emph{artifacts} in the solution of this problem. The presence of artifacts is an issue that affects many of the spectral estimation methods proposed in the literature (see, e.g., \cite[Sec.~VII-B]{FMP12}). Artifacts are usually present in the form of high and narrow frequency peaks in the spectral estimate. In light of Remark \ref{rem:art}, the advantage of using either Thompson or Hilbert distance in solving the spectral estimation problem consists of the fact that these artifacts are highly penalized by these distances, and, consequently, they should not appear in the optimal spectral estimate.

More generally, the optimization of $\H_\infty$-norms in place of or in complement to $\H_{2}$-norms has been a very fruitful direction of research in linear system theory. Building upon this heritage, the distances introduced in the present paper could also open  novel avenues in robust statistical estimation. This program raises a number of open questions, as Finsler optimization is a far less mature research area than Riemannian optimization.

Finally, Finsler distances are also a promising tool to introduce a geometry in the set of passive systems whose transfer functions are {\em positive real} and hence are in one-to-one correspondence with a spectral density.
Negative imaginary systems may inherit this geometry in view of the  connections between positive real and negative imaginary transfer functions \cite{7317760,FERRANTE20171}.







\bibliographystyle{IEEEtran}
\bibliography{Bibliography}

\end{document}